\documentclass[11pt,a4paper]{article}
\usepackage[]{times}
\usepackage{fullpage}
\usepackage{graphicx}
\usepackage{float}
\usepackage{subfigure}
\usepackage{amsmath}
\usepackage{fancyhdr}
\usepackage[colorlinks]{hyperref}
\usepackage{xcolor}
\usepackage{lineno}
\date{}%leave empty

\newtheorem{e-proposition}[theorem]{Proposition}

\newtheorem{e-definition}[theorem]{Definition\rm}
\newtheorem{remark}{\it Remark\/}

\title{Variants of Raviart-Thomas mixed finite elements \\ for curved domains using straight-edged tetrahedra}

\author{
    Vittoriano Ruas$^{1}$\thanks{Sorbonne Universit\'e, Campus Pierre et Marie Curie, 4 place jussieu, Couloir 55-65, 4\`eme \'etage, 75005 Paris, France.}
		\\[1mm]
  {\small $^{1}$ Institut Jean Le Rond d'Alembert, CNRS UMR 7190, 
  Sorbonne Universit\'e, Paris, France.}\\[1mm]
  {\small e-mail: {\it vitoriano.ruas@upmc.fr}}}

\begin{document}

\maketitle
\vspace{-6mm}
\hspace{7cm} 
\textbf{\small Dedicated to the memory of Gustavo Buscaglia}\\
\begin{abstract}
A numerical study of the tetrahedral Raviart-Thomas mixed finite element method is presented in the solution of model second order boundary value problems posed in a curved spatial domain. An emphasis is given to the case where normal fluxes are prescribed on a boundary portion.  In this case the question on the best way to enforce known boundary degrees of freedom is raised. It seems intuitive that the normal component of the flux variable should preferably not take up corresponding prescribed values at nodes shifted to the boundary of the approximating polyhedron in the underlying normal direction. This is because an accuracy downgrade is to be expected, as shown in \cite{FBRT} and \cite{ESAIM-M2AN}. While in the former work accuracy improvement is achieved by means of a standard Galerkin formulation with parametric elements, the latter one advocates the use of a straight-edged triangles combined with a Petrov-Galerkin formulation in which the test-flux space is a little different from the shape-flux space. The first purpose of this article is to show that the method studied in \cite{ESAIM-M2AN} for two-dimensional problems can be extended quite naturally to the three-dimensional case. Moreover, we carry out numerical experimentation with such a version for the two lowest order methods of this family, as compared to other strategies. In the case of the lowest order method this comparative study is enriched by assessing as well the performance of its Hermite analog 
introduced in \cite{JCAM}.


\end{abstract}

\section{Introduction}

This work deals with a specific type of Petrov-Galerkin formulation, designed to preserve orders higher than one in natural norms, and/or to improve the accuracy of finite element methods to solve boundary value problems posed in domains having a smooth curved boundary, on which degrees of freedom (DOFs) are prescribed. In previous work the author applied it to conforming or nonconforming Lagrange and Hermite finite element methods to solve both two- and three-dimensional second and fourth order elliptic equations with straight-edged triangles or tetrahedra. For all these contributions the author refers to \cite{JCAMbis} and references therein. While such a technique is simple to implement, it definitively eliminates the need for curved elements and consequently the use of non affine mappings. \\
\indent Here we examine possible advantages of this formulation, in the framework of the solution of second-order boundary value problems in a smooth spatial domain, with Neumann conditions prescribed on a boundary portion, by Raviart-Thomas mixed finite elements for $N$-simplexes. %We recall that in this case the normal components of the flux variable in the underlying mixed formulation are prescribed. 
Our premise is that the approximate normal flux should preferably not take up corresponding prescribed values at nodes shifted to the boundary of the approximating polyhedron. %Moreover, even when there is no node shift, presumably the prescribed flux on the boundary should be normal to the true boundary and not to the former. 
Actually a technique following this principle was studied in \cite{FBRT} for Raviart-Thomas mixed elements, based on a parametric version of theirs with curved simplexes. In that work the authors showed that globally their approach does bring about better accuracy. In contrast, in this paper the aforementioned Petrov-Galerkin variant with straight-edged tetrahedra is employed to tackle this problem. Our experiments advocate in favor of this approach as a rule, although in some cases no significant accuracy improvement was observed.

\section{The model problem}

\indent Our numerical study in connection with the Raviart-Thomas mixed finite elements, will be carried out for the linear reaction-diffusion first order system, with a unit diffusion coefficient and a reaction coefficient $\nu$ possibly equal to zero. In this framework we observe that accuracy enhancement is strongly dependent on the solution regularity. Hence although our technique can be applied to mixed Dirichlet-Neumann boundary conditions in very general situations, in order to make sure that the solution of the problem at hand has the required regularity for its theoretical order to prevail, we consider the following model equation. \\
Let $\Omega$ be a three-dimensional smooth domain and $\Gamma$ be its boundary, consisting of two non intersecting portions $\Gamma_0$ and $\Gamma_1$ with $meas(\Gamma_1 ) >0$. We denote by ${\bf n}$ the outer normal vector to $\Gamma$ and by $V$, either the space $L^2(\Omega)$ if $\nu >0$ or $\Gamma_0 \neq \emptyset$ if $\nu =0$, or its subspace $L^2_0(\Omega)$ consisting of those functions $g$ such that $\int_{\Omega} g =0$ otherwise. 
Now given $f \in V$ the problem to solve is,
\begin{equation}
\label{Mixed}
\left\{
\begin{array}{ll}
\mbox{Find } ({\bf p}; u) \mbox{ with } \int_{\Omega} u = 0 \mbox{ if } \Gamma_0 = \emptyset \mbox{ and } \nu=0, & \mbox{such that}   \\ 
 - \nabla \cdot {\bf p} + \nu u = f \mbox{ and } {\bf p} - \nabla u = {\bf 0}                  & \mbox{ in } \Omega; \\
 u = 0 \mbox{ on } \Gamma_0 \mbox{ and } {\bf p} \cdot {\bf n}= 0 \mbox{ on } \Gamma_1 &.
\end{array}
\right.
\end{equation} 
Referring to \cite{RaviartThomas}, let ${\bf Q}$ be the subspace of ${\bf H}(div,\Omega)$ of those fields ${\bf q}$ such that 
${\bf q} \cdot {\bf n}=0$ on $\Gamma_1$. Then denoting by $(\cdot,\cdot)$ the standard inner product of $L^2(\Omega)$, problem (\ref{Mixed}) can be recast in the following equivalent variational form: 
 \begin{equation}
\label{Varimix}
\left\{
\begin{array}{l}
\mbox{Find } ({\bf p}; u) \in {\bf Q} \times V \mbox{ such that,}\\ 
 -(\nabla \cdot {\bf p}, v) + \nu(u,v) = (f,v)\; \forall v \in V; \\
 ({\bf p},{\bf q}) + ( u, \nabla \cdot {\bf q}) = 0 \; \forall {\bf q} \in {\bf Q}.
\end{array}
\right.
\end{equation} 
%It is well known that the existence and uniqueness of a solution to \eqref{Varimix} is guaranteed.\\
Our numerical study applies to the Raviart-Thomas mixed finite element method for tetrahedra \cite{RaviartThomas} known as $RT_k$,  confined to the cases where $k=0$ and $k=1$.\\
\indent We recall that, provided $u$ belongs to $H^{k+2}(\Omega)$ (see e.g. \cite{Ciarlet}), this method is of order $k+1$ in the natural norm $||| \cdot |||$ of the space ${\bf H}(div,\Omega) \times L^2(\Omega)$. More concretely, denoting the standard norm of $L^2(\Omega)$ by $\| \cdot \|$, the norm $||| \cdot |||$ is defined by 
$||| ({\bf q};v) ||| = \left[ \| {\bf q} \|^2 + \| \nabla \cdot {\bf q} \|^2 + \| v \|^2 \right]^{1/2}$. 
From the assumption that $\Gamma_0$ and $\Gamma_1$ have no common points, the above regularity of $u$ holds if $\Gamma$ is of the $C^{k}$-class and $f \in H^{k}(\Omega)$, which will be assumed henceforth. 

\section{Method description}

Before going into the description of our numerical approach to solve \eqref{Varimix}, we give some notations related to the finite-element meshes.

\subsection{Meshes and related sets}

\hspace{4mm} Let $\{{\mathcal T}_h\}_h$ be a regular family of meshes in the sense of \cite{Ciarlet}, consisting of straight-edged tetrahedra satisfying the usual compatibility conditions for the finite element method, $h$ being the maximum edge length of all the elements in ${\mathcal T}_h$. Every tetrahedron of this mesh is considered to be a closed set and ${\mathcal T}_h$ is assumed to fit $\Omega$ in such a way that all the vertexes of the polyhedron $\Omega_h$ lie on $\Gamma$, where $\Omega_h$ is the interior of $\cup_{T \in {\mathcal T}_h} T$. %We further define $\tilde{\Omega}_h := \Omega \cup \Omega_h$ and $\Omega^{'}_h := \Omega \cap \Omega_h$. 
The boundary of $\Omega_h$ %and $\tilde{\Omega}_h$ are respectively 
is denoted by $\Gamma_h$ and ${\Gamma}_{1,h}$ is the portion of $\Gamma_h$ having a non empty intersection with $\Gamma_1$.%and moreover we set $\Gamma^{'}_h:= \Omega_h \cap \Gamma$. 
%The maximum edge length of every $T \in {\mathcal T}_h$ is denoted by $h_T$. 
We assume that any element in ${\mathcal T}_h$ has at most one face contained in $\Gamma_h$.  \\
Let ${\mathcal T}_{1,h}$ be the subset of ${\mathcal T}_h$ consisting of tetrahedra $T$ having 
one face on $\Gamma_{1,h}$, say $F_T$. Referring to Figure 1, for every $T \in {\mathcal T}_{1,h}$ we denote by $O_T$ the vertex of $T$ not belonging to $\Gamma_1$,  
%; moreover we define $T^{\Gamma}$ to be the curved triangle delimited by $\Gamma$ and the two edges of $T$ intersecting at $P_T$. 
%In view of the above assumption, the interior of every triangle in ${\mathcal T}_h \setminus {\mathcal S}_h$ has an empty  intersection with $\Gamma_h$.
%We also need some additional definitions regarding the skin $(\Omega \setminus \Omega_h) \cup (\Omega_h \setminus \Omega)$. 
and denote by ${\bf n}_T$ the unit outer normal vector to $\Gamma_{1,h}$. 
%we also define $\Delta_T$ to be the closed set delimited by $\Gamma_1$ 
and by $F_T$ the face of $T$ whose vertexes belong to $\Gamma_1$.  
%In this manner we can assert that, if $\Omega$ is convex, $\Omega_h$ is a proper subset of $\Omega$ and $\bar{\Omega}$ is the union of the disjoint sets $\Omega_h$ and $\displaystyle \cup_{T \in {\mathcal S}_h} \Delta_T$. Otherwise $\Omega_h \setminus \Omega$ is a nonempty set containing subsets of $T \in {\mathcal S}_h$ whose area is an $O(h^3)$ corresponding to non-convex portions of $\Gamma$. \\%Whatever the case, the above configurations are of  merely academic interest and carry no practical meaning, as much as the sets 
%$T^{\Delta}:=T \cup \Delta_T$ and $T^{'}:= T \cap \Omega$ $\forall T \in {\mathcal S}_h$.\\
%Notice that if $e_T$ lies on a convex portion of $\Gamma_h$, $T$ is a proper subset of $T^{\Delta}$, while $T^{'}$ is a proper subset of $T$ if $e_T$ lies in a concave portion of $\Gamma_h$. \\ 

\begin{figure}[h!]
\label{fig1}
\centerline{\includegraphics[width=3.8in]{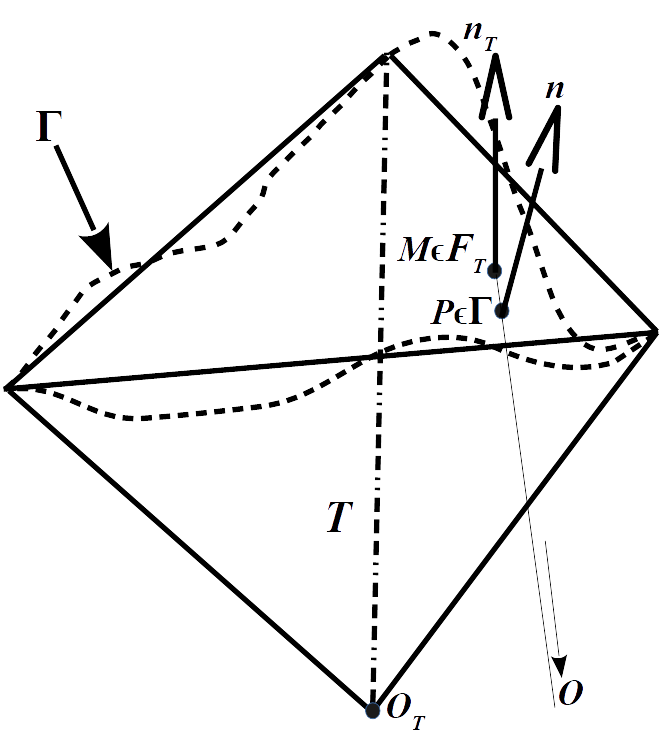}}
\vspace*{8pt}
\caption{A quadrature point $M \in F_T \subset \Gamma_{1,h}$ and normal vectors at $M$ and $P \in \overline{OM} \cap \Gamma$ for $T \in {\mathcal T}_{1,h}$}
\end{figure} 

\subsection{A Petrov-Galerkin mixed formulation}
\indent Let us introduce three spaces ${\bf P}_h^k$, ${\bf Q}_h^k$ and $V_h^k$ associated with ${\mathcal T}_h$. \\
$\bullet$ $V_h^k$ is the space consisting of functions $v \in L^2(\Omega_h)$, whose restriction to every $T \in {\mathcal T}_h$ is a polynomial of degree less than or equal to $k$ for $k \geq 0$, that vanish at a given point of $\Omega_h$ in case $\Gamma_0 =\emptyset$. \\%We extend every function $v \in V_h$ by zero in $\Omega \setminus \Omega_h$ . \\
$\bullet$ ${\bf Q}_h^k$ is the space of fields in ${\bf q} \in {\bf H}(div,\Omega_h)$, whose restriction to every $T \in {\mathcal T}_h$ belongs to the standard Raviart-Thomas space $RT_k$ \cite{RaviartThomas}, fulfilling ${\bf q} \cdot {\bf n}_T=0$ on 
$F_T$  for every $T \in {\mathcal T}_{1,h}$. \\
$\bullet$ ${\bf P}_h^k$ is the space of fields ${\bf r}$ satisfying the following conditions:\\
%With every $T \in {\mathcal S}_h$ we associate a skin $\Delta_T$ defined as follows.   It is clear that the union of all elements in ${\mathcal T}_h^{'}$ is the closure of the set $\Omega_h \cup \Omega$, and in case $\Omega$ is convex it is nothing but the closure of this set itself.
%\begin{enumerate} 
\hspace*{2mm} 1.  
${\bf r} \in {\bf H}(div,\Omega_h)$;\\
\hspace*{2mm} 2. The restriction of ${\bf r}$ to every $T \in {\mathcal T}_h$ is a field belonging to the standard Raviart-Thomas space $RT_k$ \cite{RaviartThomas} of degree less than or equal to $k$;
\\
\hspace*{2mm} 3. ${\bf r}$ is also defined in $\bar{\Omega} \setminus \bar{\Omega}_h$ in such a way that the expression of ${\bf r}$ in $T$ extends to points in every subset of the set $\Delta_T$ delimited by $F_T$, $\Gamma$ and three arbitrary plane sets whose boundary is the union of an edge of $F_T$ and a subset of $\Gamma$ for $T \in {\mathcal T}_{1,h}$ \footnote{However, akin to \cite{IMAJNA}, for any edge common to two faces in $\Gamma_{1,h}$ the two sets corresponding to each one of them must coincide.};\\ 
\hspace*{2mm} 4. $\forall T \in {\mathcal T}_{1,h}$, ${\bf r}(P) \cdot {\bf n} = 0$ for all $P$ which is the nearest intersections with $\Gamma$ of the half-line with origin at a certain point $O$ in the convex hull of $\Omega$ and $m_k:=(k+2)(k+1)/2$ points $M$ on the face $F_T$, whose minimum distance to each other is greater than $C h_T$ for a suitable mesh-independent constant $C>0$.

\begin{remark} %The construction of the nodes $P \in \Gamma_1$ associated with ${\bf P}_h^k$ advocated in item 4. is not mandatory. While it is a balanced choice in all cases addressed hereafter,   % Notice that it differs 
%from the intuitive choice of nodes lying on normals to edges of $\Gamma_{1,h}$. The main advantage of this proposal is an easy determination of boundary-node coordinates by linearity, using a supposedly available analytical expression of $\Gamma_1$. 
In principle, any way to choose such $m_k$ nodes would be all alright for our Petrov-Galerkin formulation to work. However, a balanced choice in terms of accuracy and order preservation is the set of $m_k$ quadrature points of a Gaussian formula that possibly integrates exactly polynomials defined on $F_T$ of degree less than or equal to $k+1$. For more details on such formulas we refer to \cite{Dunavant} and \cite{Stroud}. Incidentally, in our numerical studies we chose $M$ to be the centroid of $F_T$ for $k=0$ and the points with barycentric coordinates of $F_T$ equal to $(2/3,1/6,1/6)$, $(1/6,2/3,1/6)$ and $(1/6,1/6,2/3)$ for $k=1$. $O$ in turn is the centroid of $\Omega$ in all cases.  \rule{2mm}{2mm}                                                                                                                                       
\end{remark}

\indent It is possible to prove that ${\bf P}_h^k$ is a non empty finite-dimensional space, provided $h$ is sufficiently small. This result can be established following the main lines of \cite{ESAIM-M2AN} for its two-dimensional analog. \\
   
Now let us set the problem associated with spaces ${\bf P}_h^k$, ${\bf Q}_h^k$ and $V_h^k$, whose solution is an approximation of the solution of (\ref{Varimix}). First we extend $f$ to $\Omega_h \setminus \Omega$ if applicable, by requiring that such an extension belongs to $H^{k}(\Omega \cup \Omega_h)$. Still denoting the resulting function by $f$ and defining $(\cdot,\cdot)_h$ to be the standard inner product of $L^2(\Omega_h)$, the problem to solve is,
\begin{equation}
\label{Varimixh}
\left\{
\begin{array}{l}
\mbox{Find } ({\bf p}_h;u_h) \in {\bf P}_h^k \times V_h^k \mbox{ such that} \\
-(\nabla \cdot {\bf p}_h,v)_h + \nu (u_h,v)_h = (f,v)_h \; \forall v \in V_h^k; \\ 
({\bf p}_h,{\bf q})_h + (u_h,\nabla \cdot {\bf q})_h = 0 \; \forall {\bf q} \in {\bf Q}_h^k.
\end{array}
\right.
\end{equation}
  
\subsection{Petrov-Galerkin formulation for the Hermite analog of $RT_0$}

Referring to \cite{JCAM} for a Hermite version of $RT_0$, we briefly recall it below as applied to problem \eqref{Varimix}.\\
Let $P_2(D)$ denote the space of polynomials of degree less than or equal to two in a bounded set $D$ of $\Re^3$ and let $\nabla_h$ denote  the broken gradient of a function $w$ such that $\nabla w_{|T} \in [L^2(T)]^3 \; \forall T \in {\mathcal T}_h$, that is, $(\nabla_h w)_{|T}:= \nabla w_{|T} \; \forall T \in {\mathcal T}_h$. Similarly, for every function $w$ having the above property that $w$ fulfills $\Delta w_{|T} \in L^2(T) \; \forall T \in {\mathcal T}_h$ as well, we define its broken laplacian by $(\Delta_h w)_{_T} = \Delta w_{|T} \; \forall T \in {\mathcal T}_h$. Keeping in mind these notations and confining ourselves to the case where $\Gamma_0 \neq \emptyset$ for simplicity, we define:\\

\noindent $\bullet$ $V^H_h := \{w | \; w \in L^2(\Omega_h), \; w_{|T} \in P_2(T) \; \forall T \in {\mathcal T}_h \mbox{ and } \nabla_h w \in {\bf Q}_h^0\}$;\\
$\bullet$ $U^H_h := \{w | \; w \in L^2(\Omega_h), \; w_{|T} \in P_2(T) \; \forall T \in {\mathcal T}_h \mbox{ and } \nabla_h w \in {\bf P}_h^0\}$. \\

\noindent Then we approximate $u$ by $u_h^H \in U^H_h$ and ${\bf p}$ by $\nabla_h u_h^H$, where $u_h^H$ is the unique function such that:
 \begin{equation}
\label{HRT0}
\left\{
\begin{array}{l}
a_h(u_h^H,v) = (f,v)_h \; \forall v \in V_h^H, \mbox{ with}\\
a_h(w,v) := (\nabla_h w, \nabla_h v)_h - (\Delta_h w, v)_h + (w, \Delta_h v)_h + \nu(w,v)_h \; \forall w \in U^H_h \mbox{ and } \forall v \in V^H_h.
\end{array}
\right.
\end{equation} 
The main difference between \eqref{HRT0} and \eqref{Varimixh} for $k=0$ is the fact that $u$ is approximated in a space of incomplete quadratic functions in each tetrahedron, though containing all linear functions. As demonstrated in \cite{JCAM}, this formally leads to second order approximations of $u$ in $L^2(\Omega)$, whenever $\Omega$ is a convex polygon, while maintaining first order approximations of ${\bf p}$ by $\nabla_h u^H_h$ in ${\bf H}(div,\Omega)$. A similar behavior is expected in the case of smooth domains, which will be confirmed hereafter.\\
Incidentally, in the sequel we will refer to the solution method inherent to problem \eqref{HRT0} as the $HRT0$-method

\subsection{Well-posedness}
%%%%%%%%%%%%%%%%%%%%%%%%%%%%%%%%%%%%%%%%%%%%%%%%%%%%%%%%%%%%%%%%%%%%%%%%%%%%%%%%%%%%%%%%%%%%%%%%%%%%%%%%%%%%%%%%%%%%%%%%%%%%%%%%%%%
We next assume that both problems \eqref{Varimixh} and \eqref{HRT0} 
have a unique solution. Referring to \cite{ESAIM-M2AN} this assumption should turn out to hold at least as long as $h$ is sufficiently small. As a consequence, as far as \eqref{Varimixh} is concerned, the 
underlying method is expected to converge with order $k+1$ in the norm of ${\bf H}(div, \Omega_h) \times L^2(\Omega_h)$, as $h$ goes to zero, as long as the method is mesh-independently (i.e. uniformly) stable in the natural norm (cf. \cite{ESAIM-M2AN}). As for \eqref{HRT0}, under the same conditions it can be expected that the Hermite method converges with order two in the norm of $L^2(\Omega_h)$ and with order one in the broken semi-norm of $[H^1(\Omega_h)]^3$, while the approximate solution's broken Laplacian converges to $\Delta u$ with order one in $L^2(\Omega_h)$. \\
\indent In the remainder of this article we attempt to validate the above expectations for both \eqref{HRT0} and \eqref{Varimixh} for $k=0$ together with the variant \eqref{Varimixh1} of the latter for $k=1$, to be defined in Subsection 4.2. Furthermore we compare the performance of these three methods designed to deal more properly with normal flux DOFs prescribed on curvilinear boundaries, with corresponding do-nothing strategies. More specifically, the latter mean replacing ${\bf P}_h^0$ with ${\bf Q}_h^0$ in \eqref{Varimixh} and $U_h^H$ with $V_h^H$ in \eqref{HRT0}, or yet ${\bf P}_h^1$ with ${\bf Q}_h^1$ in \eqref{Varimixh1}. In the sequel the approximate problems related to such a do-nothing solution strategy will be referred to as Galerkin formulations.
%%%%%%%%%%%%%%%%%%%%%%%%%%%%%%%%%%%%%%%%%%%%%%%%%%%%%%%%%%%%%%%%%%%%%%%%%%%%%%%%%%%%%%%%%%%%%%%%%%%%%%%%%%%%%%%%%%%%%%%%%%%%%%%%%%%
   
\section{Numerical experiments}

In this section we go into numerical experimentation, which is the main purpose of this paper. At first we do this with the Petrov-Galerkin formulation defined by \eqref{Varimixh} for $k=0$ and also for the Hermite variant of the latter stipulated in \eqref{HRT0}. Then we pursue the experiments using the alternative formulation \eqref{Varimixh1} with $\nu >0$ for $k=1$. We check out the performances of the three methods by solving test-problems with known exact solution satisfying different types of boundary conditions. In all cases $\Omega$ is a domain contained in the ellipsoid $\Omega_1$ with semi-axes equal to $a$, $b$ and $c$, whose boundary is denoted by $\partial \Omega_1$. More concretely, we consider both the case where $\Omega$ is $\Omega_1$ itself or the case where $\Omega$ is the hollow domain $\Omega_{hol} := \Omega_1 \setminus \overline{\Omega_{1/2}}$, $\Omega_{1/2}$ being the concentric ellipsoid with boundary $\partial \Omega_{1/2}$ and semi-axes $a/2$, $b/2$ and $c/2$ aligned with those of $\Omega_1$. Since we only considered test-problems with three planes of symmetry, we computed with quasi-uniform family of meshes for the octant ${\Omega}_1^{+++}$ of $\Omega_1$ given by $x>0$, $y>0$ and $z>0$ generated by the procedure described in \cite{LBRSR} defined by an even integer $L$. We recall that for $\Omega_1$ such meshes consist of $6L^3$ tetrahedra. In case $\Omega$ is the hollow domain $\Omega_{hol}$ the final mesh is obtained by simply removing from the mesh of ${\Omega}_1^{+++}$ the $3L^3/4$ tetrahedra fully contained in the corresponding octant of the ellipsoid $\Omega_{1/2}$. \\
\indent Henceforth we use as mesh parameter the quantity $h:=1/L$. In all cases we take $L=2^m$ for $m$ ranging between $1$ and $4$.\\

For all the test-problems addressed in this section the results are summarized in tables with the following layout:\\
First we note that a given table presents results obtained with either a Petrov-Galerkin formulation such as \eqref{Varimixh} and \eqref{HRT0}, or their Galerkin counterparts, as indicated in the table caption. \\
Each table consists of three subsequent pairs of rows in which the approximation errors of ${\bf p}$, $\nabla \cdot {\bf p}$ and $u$ measured in the standard norm of $L^2(\Omega_h)$ for increasing values of $m$ are displayed. In the first row of each pair we display the errors for the methods introduced in this work. In order to highlight eventual advantages of our Petrov-Galerkin approach over the Galerkin one, in the second row corresponding errors are supplied for the solution strategy in which ${\bf Q}_h^k$ replaces ${\bf P}_h^k$ and $V^H_h$ replaces $U^H_h$. For the sake of clarity we denote the underlying solutions by $(\bar{\bf p}_h;\bar{u}_h)$ in the case of \eqref{Varimixh} or by $\bar{u}_h^H$ in the case of \eqref{HRT0}.\\
Finally, in the last column of each row an estimated order of convergence ($EOC$) for the corresponding quantity is displayed. The order estimation is performed by linear regression from the logarithmic sample of the four errors vs. log $h$. \\

Hereafter $\| \cdot \|_h$ denotes the norm of $L^2(\Omega_h)$ and $\eta$ stands for $\sqrt{(x/a)^2+(y/b)^2+(z/c)^2}$.   

\subsection{Solution of test-problems with the $RT_0$-method or its Hermite variant}

To begin with we clarify that two approaches were used to solve the linear system of algebraic equations leading to the solution of problems \eqref{Varimixh}  and \eqref{HRT0}, according to the case, namely, \\
\begin{itemize} 
\item If $\nu =0$ Crout's method with partial pivoting for the whole system was employed; 
\item If $\nu \neq 0$ elimination by static condensation of the DOFs internal to a tetrahedron - i.e. the mean values of $u$ therein - is first carried out. The resulting linear system for the flux DOFs is solved by Crout's method with partial pivoting and afterward the internal DOFs are recovered from the former element by element.
\end{itemize}

\begin{enumerate}

\item
\underline{\bf Test-problem 1:} $\Omega=\Omega_1$; $a=b=c=1$; $u=(\eta^2/2-\eta^3/3)$; $\nu=1$; $\Gamma_1=\Gamma$. 

\indent Table 1 indicates that, akin to its Galerkin counterpart, the Petrov-Galerkin formulation \eqref{Varimixh} for $k=0$ yields as expected first order approximations of the flux field ${\bf p}$ in the norm of ${\bf H}(div,\Omega_h)$, which together with $u$ solves the coercive problem \eqref{Varimix} with $\nu >0$. On the other hand, surprisingly enough, both approximations $u_h$ and $\bar{u}_h$ of $u$ appear to be superconvergent in the norm of $L^2(\Omega_h)$, respectively of orders equal to ca. $1.48$ and $1.58$. Taking into account that all the computed approximations by \eqref{Varimixh} of $u$ are fairly more accurate than those obtained with the standard Galerkin formulation, and moreover the approximation errors of ${\bf p}$ in the ${\bf H}(div)$-norm by both methods are practically identical, we can claim the superiority of our Petrov-Galerkin formulation over its Galerkin counterpart, as far as the current test-problem is concerned.\\
            
\begin{table*}[h!]
{\small 
\centering
\begin{tabular}{ccccccc} &\\ [-.3cm]  
$h$ & $\longrightarrow$ & $1/2$ & $1/4$ & $1/8$ & $1/16$ & $EOC$ 
\tabularnewline &\\ [-.3cm] \hline &\\ [-.3cm]
$\parallel {\bf p}_h - {\bf p} \parallel_h$ & $\longrightarrow$ & 0.62009878E-1 & 0.30260267E-1 & 0.15024034E-1 & 0.74926747E-2 & 1.0157 
\tabularnewline &\\ [-.3cm] \hline &\\ [-.3cm] 
$\parallel \bar{\bf p}_h - {\bf p} \parallel_h$ & $\longrightarrow$ & 0.63726346E-1 & 0.30633971E-1 & 0.15079949E-1 & 0.75001745E-2 & 1.0283
\tabularnewline &\\ [-.3cm] \hline &\\ [-.3cm] 
\tabularnewline &\\ [-.3cm] \hline &\\ [-.3cm]
$\parallel \nabla \cdot({\bf p}_h - {\bf p}) \parallel_h$ & $\longrightarrow$ & 0.28540617E+0 & 0.14667675E+0 & 0.73739856E-1 & 0.36913166E-1 & 0.9845
\tabularnewline &\\ [-.3cm] \hline &\\ [-.3cm]  
$\parallel \nabla \cdot(\bar{\bf p}_h - {\bf p}) \parallel_h$ & $\longrightarrow$ & 0.28687043E+0 & 0.14695449E+0 & 0.73779617E-01 & 
0.36918413E-1 & 0.9868
\tabularnewline &\\ [-.3cm] \hline &\\ [-.3cm]
\tabularnewline &\\ [-.3cm] \hline &\\ [-.3cm] 
 $\parallel u_h - u \parallel_h$ & $\longrightarrow$ & 0.35280503E-1 & 0.93274239E-2 & 0.34817973E-2 & 0.16016971E-2 & 1.4805 
\tabularnewline &\\ [-.3cm] \hline &\\ [-.3cm] 
$\parallel \bar{u}_h - u \parallel_h$  & $\longrightarrow$ & 0.45636307E-1 & 0.12982795E-1 & 0.42412857E-2 & 0.17183781E-2 &  1.5807
\tabularnewline &\\ [-.3cm] \hline &\\ [-.3cm] 
\end{tabular}
%\vskip.2cm 
\caption{Errors for Test-problem 1 solved by the $RT_0$-method} 
}
\label{table1}
\end{table*}

\indent Now watching Table 2 we conclude that the situation is even more in favor of our Petrov-Galerkin formulation, for 
superconvergence of the approximations of ${\bf p}$ of order ca. $1.24$ is observed, while its standard Galerkin counterpart remains roughly equal to $1$. On the other hand the observed orders of approximation of both $u_h^H$ and $\bar{u}_h^H$ in the $L^2$-norm are roughly equal to $2$ as expected, although the accuracy of the former is significantly better all the way, as compared to the latter. It is also  noticeable that the  approximations of $\nabla \cdot {\bf p}$ are practically identical for all the four methods assessed in this test-problem. \\

In short, on the basis of the above comments regarding the current test-problem, we can assert that the $HRT_0$-method is by far the best option to solve reaction-diffusion equations with Neumann boundary conditions prescribed on the boundary of a curved spatial domain. 
   
\begin{table*}[h!]
{\small 
\centering
\begin{tabular}{ccccccc} &\\ [-.3cm]  
$h$ & $\longrightarrow$ & $1/2$ & $1/4$ & $1/8$ & $1/16$ & $EOC$ 
\tabularnewline &\\ [-.3cm] \hline &\\ [-.3cm]
$\parallel \nabla_h(u^H_h - u) \parallel_h$       & $\longrightarrow$ & 0.60342191E-1 & 0.25329447E-1 & 0.10534240E-1 & 0.46280183E-2 & 1.2380
\tabularnewline &\\ [-.3cm] \hline &\\ [-.3cm] 
$\parallel \nabla_h(\bar{u}^H_h - u) \parallel_h$ & $\longrightarrow$ & 0.68031204E-1 & 0.34550245E-1 & 0.16747424E-1 & 0.81488648E-2 & 1.0229
\tabularnewline &\\ [-.3cm] \hline &\\ [-.3cm] 
\tabularnewline &\\ [-.3cm] \hline &\\ [-.3cm]
$\parallel \Delta_h(u^H_h - u) \parallel_h$       & $\longrightarrow$ & 0.28544318E+0 & 0.14617124E+0 & 0.73740806E-1 & 0.36913286E-1 & 0.9845
\tabularnewline &\\ [-.3cm] \hline &\\ [-.3cm]  
$\parallel \Delta_h(\bar{u}^H_h - u) \parallel_h$ & $\longrightarrow$ & 0.28692715E+0 & 0.14696374E+0 & 0.73780839E-1 & 0.36918568E-1 & 0.9869
\tabularnewline &\\ [-.3cm] \hline &\\ [-.3cm]
\tabularnewline &\\ [-.3cm] \hline &\\ [-.3cm] 
$\parallel u^H_h - u \parallel_h$                & $\longrightarrow$ & 0.33901594E-1 & 0.72607746E-2 & 0.16389278E-2 & 0.38660695E-3 & 2.1510 
\tabularnewline &\\ [-.3cm] \hline &\\ [-.3cm] 
$\parallel \bar{u}^H_h - u \parallel_h$          & $\longrightarrow$ & 0.44737867E-1 & 0.11618171E-1 & 0.29314026E-2 & 0.73449907E-3 & 1.9772
\tabularnewline &\\ [-.3cm] \hline &\\ [-.3cm] 
\end{tabular}
%\vskip.2cm 
\caption{Errors for Test-problem 1 solved by the $HRT_0$-method and its Galerkin counterpart} 
}
\label{table1bis}
\end{table*}

\item
\underline{\bf Test-problem 2:} $\Omega=\Omega_1$; $a=b=c=1$; $u=(\eta^2/2-\eta^3/3)$; $\nu=0$; $\Gamma_1=\Gamma$.\\

\indent First we recall that we avoided enforcing the unpractical condition $u_h \in L^2_0(\Omega_h)$ (resp. $\bar{u}_h \in L^2_0(\Omega_h)$)in order to fix the additive constant up to which $u_h$ (resp. $\bar{u}_h$) is defined. Instead we fixed this constant by prescribing $u_h(T_0)=0$ (resp. $\bar{u}_h(T_0)=0$), where $T_0$ is one out of the six tetrahedra having the origin as a vertex.\\
\indent From Table 3 we infer that in the case under study our Petrov-Galerkin approach is more accurate than the standard Galerkin one in all respects. It is even possible to observe a significant superconvergence of the approximation of $u$ by the former method, while the expected first order convergence to ${\bf p}$ is confirmed for both methods. On the other hand, there is a severe erosion of the EOCs, as far as the approximations of $\nabla \cdot {\bf p}$ are concerned, even though to a much lesser extent for our Petrov-Galerkin approach. As we saw, in the preceding example with $\nu >0$ such a downgrade does not occur. It could be that better approximations of the Poisson equation with Neumann conditions prescribed everywhere on the boundary would be obtained, if the aforementioned arbitrary additive constant were fixed by a scalar integral condition on $u_h$ (resp. $\bar{u}_h$), instead of fixing its value at a point (i.e., in a given tetrahedron in the case of $RT_0$). However, the enforcement of such a condition increases significantly the computational effort in a direct solution of the algebraic linear system. %Such a conjecture is corroborated by the coherent results for the following example, in which mixed boundary conditions hold. 
\begin{remark} 
Scalar integral conditions on $u_h$ and $\bar{u}_h$ could be implicitly taken into account at no extra cost, in case a suitable iterative solution procedure was employed. One such a technique is currently under study, as applied to problem \eqref{Varimixh} and its Galerkin counterpart, for $k=0$ and $\nu=0$. \rule{2mm}{2mm}
\end{remark}
Notice that the arguments given in \cite{BrezziFortin} indicate that, whatever the way the arbitrariness of the additive constant of $\bar{u}_h$ is resolved, the standard Galerkin approach leads to first order approximations in $L^2(\Omega_h)$ of all the three fields ${\bf p}$, $\nabla \cdot {\bf p}$ and $u$. That is why quality improvement of the EOC to $\nabla \cdot {\bf p}$ can be expected for both methods, as meshes finer than those used in our computations come into play. %Anyhow the order erosion observed in this test-problem is probably due to the way we fixed the arbitrary additive constant of the approximations of $u$. \\
Such a conjecture is corroborated by the fact that the two-dimensional analogs of both methods assessed in this subsection are known to be uniformly stable and of the first order in the natural norms in the solution of the Poisson equation in a curved domain with different  types of boundary conditions. These include Neumann conditions prescribed everywhere on the boundary (cf. \cite{BrezziFortin} and \cite{ESAIM-M2AN}).

\begin{table*}[h!]
{\small 
\centering
\begin{tabular}{ccccccc} &\\ [-.3cm]  
$h$ & $\longrightarrow$ & $1/2$ & $1/4$ & $1/8$ & $1/16$ & $EOC$ 
\tabularnewline &\\ [-.3cm] \hline &\\ [-.3cm]
$\parallel {\bf p}_h - {\bf p} \parallel_h$ & $\longrightarrow$ & 0.69115375E-1 & 0.31925817E-1 & 0.15401040E-1 & 0.75807731E-2 & 1.0617 
\tabularnewline &\\ [-.3cm] \hline &\\ [-.3cm] 
$\parallel \bar{\bf p}_h - {\bf p} \parallel_h$ & $\longrightarrow$ & 0.74852127E-1 & 0.34901082E-1 & 0.16351475E-1 & 0.78444990E-2 & 1.0857
\tabularnewline &\\ [-.3cm] \hline &\\ [-.3cm] 
\tabularnewline &\\ [-.3cm] \hline &\\ [-.3cm]
$\parallel \nabla \cdot({\bf p}_h - {\bf p}) \parallel_h$ & $\longrightarrow$ & 0.38377087E+0 & 0.21254516E+0 & 0.12192545E+0 & 0.74154103E-1 & 0.7917
\tabularnewline &\\ [-.3cm] \hline &\\ [-.3cm]  
$\parallel \nabla \cdot(\bar{\bf p}_h - {\bf p}) \parallel_h$ & $\longrightarrow$ & 0.43740878E+0 & 0.29027154E+0 & 0.19446875E+0 & 
0.13296500E+0 & 0.5732
\tabularnewline &\\ [-.3cm] \hline &\\ [-.3cm]
\tabularnewline &\\ [-.3cm] \hline &\\ [-.3cm] 
 $\parallel u_h - u \parallel_h$ & $\longrightarrow$ & 0.60342191E-1 & 0.25329447E-1 & 0.10534240E-1 & 0.46280183E-2 & 1.2380 
\tabularnewline &\\ [-.3cm] \hline &\\ [-.3cm] 
$\parallel \bar{u}_h - u \parallel_h$  & $\longrightarrow$ & 0.68031204E-1 & 0.34550245E-1 & 0.16747424E-1 & 0.81488648E-2 &  1.0229
\tabularnewline &\\ [-.3cm] \hline &\\ [-.3cm] 
\end{tabular}
%\vskip.2cm 
\caption{Errors for Test-problem 2 solved by the $RT_0$-method in formulations \eqref{Varimixh} and Galerkin's} 
}
\label{table2}
\end{table*}

\item
\underline{\bf Test-problem 3:} $\Omega=\Omega_{hol}$; $a=0.6$, $b=0.8$, $c=1$; $u=(2 \eta-1)(\eta-1)^2/2$; $\nu=0$; $\Gamma_1= \partial \Omega_1$.\\

Akin to Test-problem 1, this test-problem is solved by both the $RT_0-$ and the $HRT_0-$method. Here we are facing similar conditions to the case of the former, since we are also solving the Poisson equation. However now Neumann conditions are prescribed only on the outer boundary $\partial \Omega_1$ of the hollow domain $\Omega_{hol}$, while Dirichlet conditions hold on 
its inner boundary $\partial \Omega_{1/2}$. Therefore we do not need to enforce any particular condition on $u$. Certainly for this reason,  optimal convergence is observed on Tables 4 and 5 for the approximations of all the three quantities being assessed, at least at expected orders. Notice that the Petrov-Galerkin approximation $u_h$ is just a little better than the approximation $\bar{u}_h$ computed with the $RT_0$-method in standard Galerkin formulation. Moreover, surprisingly, both methods behave almost like second order ones, with $ECO$s close to $1.9$. Since there is no particular reason for this, we conjecture that such an effect is only observable for coarse meshes such as those used in our computations; these $ECO$s should gradually decrease as finer meshes are used in the estimations, but a final order strictly greater than one could apply here, akin to Test-problem 1. \\
Finally we note that all the four methods experimented in this test-problem behave quite in the same manner, as far as the flux variable ${\bf p}$ is concerned, while its computed divergence is exactly the same in all cases. These findings are somewhat disappointing in the case of the former field, since our Petrov-Glaerkin formulation seems not to bring about any significant improvement here, in contrast to the two previous test-problems. Furthermore, we have no explanation for the invariance of the approximations of $\nabla \cdot {\bf p}$. On the other hand, the $HRT_0$-method is superior to the $RT_0$-method all the way, as far as the approximations of $u$ are concerned, with an advantage of the Petrov-Galerkin formulation over the standard Galerkin one.       

\begin{table*}[h!]
{\small 
\centering
\begin{tabular}{ccccccc} &\\ [-.3cm]  
$h$ & $\longrightarrow$ & $1/2$ & $1/4$ & $1/8$ & $1/16$ & $EOC$ 
\tabularnewline &\\ [-.3cm] \hline &\\ [-.3cm]
$\parallel {\bf p}_h - {\bf p} \parallel_h$     & $\longrightarrow$ & 0.82727724E-1 & 0.34134216E-1 & 0.16334542E-1 & 0.80930093E-2 & 1.1124
\tabularnewline &\\ [-.3cm] \hline &\\ [-.3cm] 
$\parallel \bar{\bf p}_h - {\bf p} \parallel_h$ & $\longrightarrow$ & 0.81765900E-1 & 0.35042560E-1 & 0.16520939E-1 & 0.81205733E-2 & 1.1080
\tabularnewline &\\ [-.3cm] \hline &\\ [-.3cm] 
\tabularnewline &\\ [-.3cm] \hline &\\ [-.3cm]
$\parallel \nabla \cdot({\bf p}_h - {\bf p}) \parallel_h$ & $\longrightarrow$ & 0.44522265E+0 & 0.25290810E+0 & 0.12997231E+0 & 0.65419324E-1 & 0.9261
\tabularnewline &\\ [-.3cm] \hline &\\ [-.3cm]  
$\parallel \nabla \cdot(\bar{\bf p}_h - {\bf p}) \parallel_h$ & $\longrightarrow$ & 0.44522265E+0 & 0.25290810E+0 & 0.12997231E+0 & 
0.65419324E-1 & 0.9261
\tabularnewline &\\ [-.3cm] \hline &\\ [-.3cm]
\tabularnewline &\\ [-.3cm] \hline &\\ [-.3cm] 
 $\parallel u_h - u \parallel_h$ & $\longrightarrow$ & 0.35222621E-1 & 0.75289989E-2 & 0.19185813E-2 & 0.61320281E-3 & 1.9504 
\tabularnewline &\\ [-.3cm] \hline &\\ [-.3cm] 
$\parallel \bar{u}_h - u \parallel_h$  & $\longrightarrow$ & 0.34731132E-1 & 0.81310202E-2 & 0.21248225E-2 & 0.65983688E-3 & 1.9090 
\tabularnewline &\\ [-.3cm] \hline &\\ [-.3cm] 
\end{tabular}
%\vskip.2cm 
\caption{Errors for Test-problem 3 solved by the $RT_0$-method in formulations \eqref{Varimixh} and Galerkin's} 
}
\label{table3}
\end{table*}

\begin{table*}[h!]
{\small 
\centering
\begin{tabular}{ccccccc} &\\ [-.3cm]  
$h$ & $\longrightarrow$ & $1/2$ & $1/4$ & $1/8$ & $1/16$ & $EOC$ 
\tabularnewline &\\ [-.3cm] \hline &\\ [-.3cm]
$\parallel \nabla_h(u^H_h - u) \parallel_h$       & $\longrightarrow$ & 0.83726242E-1 & 0.34403190E-1 & 0.16374151E-1 & 0.80981944E-2 & 1.1181
\tabularnewline &\\ [-.3cm] \hline &\\ [-.3cm] 
$\parallel \nabla_h(\bar{u}^H_h - u) \parallel_h$ & $\longrightarrow$ & 0.82702529E-1 & 0.35300610E-1 & 0.16560026E-1 & 0.81257431E-2 & 1.1150
\tabularnewline &\\ [-.3cm] \hline &\\ [-.3cm] 
\tabularnewline &\\ [-.3cm] \hline &\\ [-.3cm]
$\parallel \Delta_h (u_h^H - u) \parallel_h$      & $\longrightarrow$ & 0.44522265E+0 & 0.25290810E+0 & 0.12997231E+0 & 0.65419324E-1 & 0.9261
\tabularnewline &\\ [-.3cm] \hline &\\ [-.3cm]  
$\parallel \Delta_h (\bar{u}^H_h - u) \parallel_h$& $\longrightarrow$ & 0.44522265E+0 & 0.25290810E+0 & 0.12997231E+0 & 0.65419324E-1 & 0.9261
\tabularnewline &\\ [-.3cm] \hline &\\ [-.3cm]
\tabularnewline &\\ [-.3cm] \hline &\\ [-.3cm] 
 $\parallel u^H_h - u \parallel_h$                & $\longrightarrow$ & 0.30621669E-1 & 0.62949113E-2 & 0.14478688E-2 & 0.34963476E-3 & 2.1478
\tabularnewline &\\ [-.3cm] \hline &\\ [-.3cm]
$\parallel \bar{u}^H_h - u \parallel_h$           & $\longrightarrow$ & 0.30088329E-1 & 0.68875925E-2 & 0.16647189E-2 & 0.41216647E-3 & 2.0618
\tabularnewline &\\ [-.3cm] \hline &\\ [-.3cm] 
\end{tabular}
%\vskip.2cm 
\caption{Errors for Test-problem 3 solved by the $HRT_0$-method and its Galerkin counterpart} 
}
\label{table3bis}
\end{table*}

\end{enumerate}

\subsection{Solution of test-problems with the $RT_1$-method}

In the three-dimensional case the mere assessment of the $RT_1$-method in solving a model problem like \eqref{Varimixh} engenders high computational costs as meshes are refined. This is mainly due to the resulting increasingly large   
non symmetric systems of algebraic equations, even at moderate levels of refinement. Thus a word is in order on the strategy that we adopted to solve them, with the aim of reducing both storage requirements and running time as much as possible. \\
Actually, such a strategy applies only to the case where $\nu >0$. For this reason in no test-problem addressed in this subsection we took $\nu =0$.\\
First of all, since $\nu \neq 0$, we may consider the following enriched form of \eqref{Varimixh}, namely, the formulation
\begin{equation}
\label{Varimixh1}
\left\{
\begin{array}{l}
\mbox{Find } ({\bf p}_h;u_h) \in {\bf P}_h^1 \times V_h^1 \mbox{ such that} \\
-(\nabla \cdot {\bf p}_h,v)_h + \nu (u_h,v)_h = (f,v)_h \; \forall v \in V_h^1; \\ 
({\bf p}_h,{\bf q})_h - (u_h,\nabla \cdot {\bf q})_h + 2 \nu^{-1} (\nabla \cdot {\bf p}_h,\nabla \cdot {\bf q})_h = -2 \nu^{-1} 
(f,\nabla \cdot {\bf q})_h \; \forall {\bf q} \in {\bf Q}_h^1.
\end{array}
\right.
\end{equation}
We recall that the Galerkin counterpart of \eqref{Varimixh1} is the problem resulting from the replacement of ${\bf P}_h^1$ with  
${\bf Q}_h^1$, whose solution is still denoted by $(\bar{\bf p}_h,\bar{u}_h)$. It is not difficult to establish that the bilinear form associated with \eqref{Varimixh1} is symmetric and coercive on $[{\bf H}(div,\Omega_h) \times L^2(\Omega_h)] \times [{\bf H}(div,\Omega_h) \times L^2(\Omega_h)]$ equipped with the natural norm. Therefore such a problem leads to symmetric linear system of algebraic equations with a positive definite matrix, for whose solution it seems judicious to apply the iterative procedure advocated in \cite{JCAMbis}% to solve \eqref{Varimixh1}. 
In the particular case under study it functions as follows:\\
Starting from a suitable pair $({\bf p}_h^0,u_h^0) \in {\bf P}_h^1 \times V_h^1$, for $n=1,2,\ldots$, determine approximations $({\bf p}_h^n,u_h^n) \in {\bf P}_h^1 \times V_h^1$ of the solution $({\bf p}_h,u_h)$ of \eqref{Varimixh1}, until a given tolerance $tol$ is satisfied in a sense to be specified, say, for $n=N_{it}$. Here the initial approximation $({\bf p}_h^0,u_h^0)$ is chosen to be the pair defined as follows:\\    
$u^0_h:=\bar{u}_h$ and ${\bf p}_h^0$ is a certain ${\bf P}_h^1$-interpolate of $\bar{\bf p}^0_h:=\bar{\bf p}_h$, where the pair $(\bar{\bf p}_h,\bar{u}_h) \in {\bf Q}_h^1 \times V^1_h$ is the solution of the Galerkin counterpart of \eqref{Varimixh1}. More specifically, referring to \cite{JCAMbis}, the pair $({\bf p}_h^n,u_h^n)$ is such that ${\bf p}_h^n$ is the ${\bf P}_h^1$-interpolate in the same sense of $\bar{\bf p}_h^n \in {\bf Q}^1_h$, where $(\bar{\bf p}_h^n,u_h^n) \in {\bf Q}^1_h \times V^1_h$ is the solution of the Galerkin counterpart of \eqref{Varimixh1} with right hand sides modified in the following manner:\\
Functional $F_h^n(v)$ replaces $(f,v)_h$ in the first equation and functional $G_h^n({\bf q})$ replaces zero in the second equation, where
$$F_h^n(v):= (f,v)_h -(\nabla \cdot [\bar{\bf p}^{n-1}_h-{\bf p}^{n-1}_h],v)_h \; \forall v \in V_h^1$$
and
$$G_h^n({\bf q}):= (\bar{\bf p}^{n-1}_h-{\bf p}^{n-1}_h,{\bf q})_h + 2 \nu^{-1}(\nabla \cdot [\bar{\bf p}^{n-1}_h-{\bf p}^{n-1}_h],\nabla \cdot {\bf q})_h \; \forall {\bf q} \in {\bf Q}_h^1.$$
\noindent The aforementioned ${\bf P}^1_h$-interpolate of a given $\bar{\bf r} \in {\bf Q}^1_h$ is the field ${\bf r} \in {\bf P}^1_h$ whose DOFs coincide with those of $\bar{\bf r}$, except for the zero values prescribed at the three quadrature points $M$ on each face contained in $\Gamma_{1,h}$, which are enforced instead at the corresponding points $P \in \Gamma$ (cf. Figure 1).\\
Now the symmetric coercive problem to determine $(\bar{\bf p}_h^n,u^n_h)$ for every $n$ can be solved by the conjugate gradient method, taking $(\bar{\bf p}_h^{n-1},u^{n-1}_h)$ as initial guess. Since the matrix associated with the problem at hand is very sparse, this amounts for a very significant reduction of storage requirements as compared to those of a direct method such as Crout's. As a matter of fact, we endeavored to reduce even more the size of the linear system of algebraic equations to solve at each step of the above iterative procedure, by using static condensation of the seven unknowns attached to the interior of each tetrahedron $T \in {\mathcal T}_h$. These unknowns are the values of $u_h^n$ at the four vertexes of $T$, together with the three moments in $T$ of order zero of $\bar{\bf p}^n_h$. In doing so, the final linear system to solve also has a positive definite matrix and as unknowns only the three normal flux DOFs per face of the mesh not contained in $\Gamma_{1,h}$. We observe that at each step and for each tetrahedron the aforementioned seven inner DOFs can be determined on an element-by-element basis, using the 133 element matrix' coefficients relating these DOFs to each other and also to the normal flux DOFs across the faces of the corresponding tetrahedron. Notice also that storing these coefficients for a given mesh can be substantially more costly than storing the non zero coefficients of the final system matrix itself. Nevertheless the total storage requirements 
for the adopted linear system solution strategy remain quite moderate in the end.\\
As for the resulting running time, we can report very drastic reductions as compared to the solution by a direct method, such as Crout's. Our approach turned out to be much less costly as well than the one based on the aforementioned step-by-step procedure, in which the underlying symmetric linear systems with a fixed positive definite matrix are solved by Cholesky's method.\\
In the tables that follow, besides the approximation errors in the $L^2$-norm of the quantities ${\bf p}_h$, $\nabla \cdot {\bf p}_h$ and $u_h$ together with those of their Galerkin counterparts, we also display the corresponding $EOC$s. Referring to the case of $RT_0$, here again these estimations are derived from the approximations obtained for $h = 1/L$ with $L=2^m$ for $m=1,2,3,4$, by linear regression.\\
Additionally, here the last row of all tables show the number of iterations $N_{it}$ necessary to attain a maximum increment of degrees of freedom of the normal flux DOFs between ${\bf p}_h^{n-1}$ and ${\bf p}^n_h$ less than or equal to $10^{-5}$. In case this tolerance is not satisfied after $100$ iterations, the value $N_{it} = 100$ is displayed.    

\begin{enumerate}

\item
\underline{\bf Test-problem 4:} $\Omega=\Omega_1$; $a=b=c=1$; $u=(\eta^2/2-\eta^3/3)$; $\nu=1$; $\Gamma_1=\Gamma$. 

\indent Table 6 indicates that the Petrov-Galerkin formulation \eqref{Varimixh1} yields approximations of both the field ${\bf p}$ in the norm of ${\bf H}(div,\Omega_h)$ and $u$ in the norm of $L^2(\Omega_h)$ of orders of about $6$ to $7\%$ less than the expected second order. On the other hand, the approximations $\bar{\bf p}_h$ of ${\bf p}$ and $\bar{u}_h$ of $u$ obtained by the standard Galerkin formulation appear to be of order closer to two in their respective natural norms. However all the computed approximations by \eqref{Varimixh1} are significantly more accurate than those obtained with its standard Galerkin counterpart, especially $\bar{u}_h$. In view of this we can claim the superiority of our Petrov-Galerkin formulation over the latter approach, as far as the current test-problem is concerned.\\
We should also note that the iterative procedure described at the beginning of this subsection behaved quite well in the present case, as one can infer from Table 6 as well. 
            
\begin{table*}[h!]
{\small 
\centering
\begin{tabular}{ccccccc} &\\ [-.3cm]  
$h$ & $\longrightarrow$ & $1/2$ & $1/4$ & $1/8$ & $1/16$ & $EOC$  
\tabularnewline &\\ [-.3cm] \hline &\\ [-.3cm]
$\parallel {\bf p}_h - {\bf p} \parallel_h$ & $\longrightarrow$ & 0.28611579E-2 & 0.79826410E-3 & 0.20646241E-3 & 0.58530915E-4 & 1.8785 
\tabularnewline &\\ [-.3cm] \hline &\\ [-.3cm] 
$\parallel \bar{\bf p}_h - {\bf p} \parallel_h$ & $\longrightarrow$ & 0.11191766E-1 & 0.29597249E-2 & 0.75022471E-3 & 0.19009492E-3 & 1.9619 
\tabularnewline &\\ [-.3cm] \hline &\\ [-.3cm] 
\tabularnewline &\\ [-.3cm] \hline &\\ [-.3cm]
$\parallel \nabla \cdot({\bf p}_h - {\bf p}) \parallel_h$ & $\longrightarrow$ & 0.15169942E-1 & 0.43969701E-2 & 0.11797263E-2 & 0.30528238E-3 & 1.8803 
\tabularnewline &\\ [-.3cm] \hline &\\ [-.3cm]  
$\parallel \nabla \cdot(\bar{\bf p}_h - {\bf p}) \parallel_h$ & $\longrightarrow$ & 0.46405212E-1 & 0.12210449E-1 & 0.31077387E-2 & 
0.78249183E-3 & 1.9644
\tabularnewline &\\ [-.3cm] \hline &\\ [-.3cm]
\tabularnewline &\\ [-.3cm] \hline &\\ [-.3cm] 
 $\parallel u_h - u \parallel_h$ & $\longrightarrow$ & 0.20358578E-2 & 0.64933204E-3 & 0.17145397E-3 & 0.43711080E-4 & 1.8546 
\tabularnewline &\\ [-.3cm] \hline &\\ [-.3cm] 
$\parallel \bar{u}_h - u \parallel_h$  & $\longrightarrow$ & 0.43902861E-1 & 0.11409793E-1 & 0.28802225E-2 & 0.72180808E-3 &  1.9766
\tabularnewline &\\ [-.3cm] \hline &\\ [-.3cm]
\tabularnewline &\\ [-.3cm] \hline &\\ [-.3cm]
$N_{it}$ for $tol = 10^{-5}$ & $\longrightarrow$ & 8 & 6 & 7 & 4 & 
\tabularnewline &\\ [-.3cm] \hline &\\ [-.3 cm] 
\end{tabular}
%\vskip.2cm 
\caption{Errors for $RT_1$-solutions of Test-problem 4 \& number of iterations to attain the solution of \eqref{Varimixh1}} 
}
\label{table4}
\end{table*}

\item
\underline{\bf Test-problem 5:} $\Omega=\Omega_1$; $a=0.6$, $b=0.8$, $c=1$; $u=(\eta^2/2-\eta^3/3)$; $\nu=1$; $\Gamma_1=\Gamma$.\\

\indent  This test-problem is similar to Test-problem 4, but now the case of an ellipsoidal domain quite different from a unit ball is assessed instead. Nonetheless the results shown in Table 7 leads to conclusions quite similar to those reported for the previous test-problem, except for two observations. On the one hand the $EOC$ of ${\bf p}_h$ is roughly as close to two as the $EOC$ of $\bar{\bf p}_h$, 
which is rather good news. On the other hand, for the finest mesh only a tolerance of ca. $10^{-4}$ could be satisfied within $100$ iterations of the procedure employed to determine the former field. But this is certainly not so serious, since not only the order of magnitude of the  DOFs of ${\bf p}$ is still much larger than such a tolerance, but also because the final computed value of ${\bf p}_h$ for $h=1/16$ turns out to be very coherent with the other ones. 

\begin{table*}[h!]
{\small 
\centering
\begin{tabular}{ccccccc} &\\ [-.3cm]  
$h$ & $\longrightarrow$ & $1/2$ & $1/4$ & $1/8$ & $1/16$ & $EOC$ 
\tabularnewline &\\ [-.3cm] \hline &\\ [-.3cm]
$\parallel {\bf p}_h - {\bf p} \parallel_h$ & $\longrightarrow$ & 0.57066172E-2 & 0.15069647E-2 & 0.38413471E-3 & 0.10014124E-3 & 1.9470 
\tabularnewline &\\ [-.3cm] \hline &\\ [-.3cm] 
$\parallel \bar{\bf p}_h - {\bf p} \parallel_h$ & $\longrightarrow$ & 0.13283767E-1 & 0.36036893E-2 & 0.91998424E-3 & 0.23310436E-3 & 1.9467
\tabularnewline &\\ [-.3cm] \hline &\\ [-.3cm] 
\tabularnewline &\\ [-.3cm] \hline &\\ [-.3cm]
$\parallel \nabla \cdot({\bf p}_h - {\bf p}) \parallel_h$ & $\longrightarrow$ & 0.25490647E-1 & 0.76902311E-2 & 0.20952956E-2 & 0.54578314E-3 & 1.8512
\tabularnewline &\\ [-.3cm] \hline &\\ [-.3cm]  
$\parallel \nabla \cdot(\bar{\bf p}_h - {\bf p}) \parallel_h$ & $\longrightarrow$ & 0.62821731E-1 & 0.16989787E-1 & 0.43743436E-2 & 
 0.11072157E-2 & 1.9436
\tabularnewline &\\ [-.3cm] \hline &\\ [-.3cm]
\tabularnewline &\\ [-.3cm] \hline &\\ [-.3cm] 
 $\parallel u_h - u \parallel_h$ & $\longrightarrow$ & 0.14256441E-2 & 0.46151533E-3 & 0.12223178E-3 & 0.31092210E-4 & 1.8474 
\tabularnewline &\\ [-.3cm] \hline &\\ [-.3cm] 
$\parallel \bar{u}_h - u \parallel_h$  & $\longrightarrow$ & 0.57435439E-1 & 0.15156724E-1 & 0.38418171E-2 & 0.96385388E-3 &  1.9671
\tabularnewline &\\ [-.3cm] \hline &\\ [-.3cm] 
\tabularnewline &\\ [-.3cm] \hline &\\ [-.3cm]
$N_{it}$ for $tol = 10^{-5}$ & $\longrightarrow$ & 10 & 7 & 23 & 100 & 
\tabularnewline &\\ [-.3cm] \hline &\\ [-.3 cm] 
\end{tabular}
%\vskip.2cm 
\caption{Errors for $RT_1$-solutions of Test-problem 5 \& number of iterations to attain the solution of \eqref{Varimixh1}}
}
\label{table5}
\end{table*}

\item
\underline{\bf Test-problem 6:} $\Omega=\Omega_{hol}$; $a=b=c=1$; $u=(2 \eta-1)(\eta-1)^2/2$; $\nu=1$; $\Gamma_1= \partial \Omega_1$.\\

As illustrated in Table 8, in this test-problem the approximations of the three quantities being assessed behave a little differently from those of the previous two test-problems. First of all there are less discrepancies between the approximations of $u$ by $u_h$ and $\bar{u}_h$, although the errors of the latter are still about four times smaller. On the other hand the quality of the approximations of ${\bf p}$ by ${\bf p}_h$ is considerably better than the one of $\bar{\bf p}_h$, as compared to the case of the two preceding test-problems. Some differences can also be noted, as far as the $ECO$s are concerned. First of all, for the standard Galerkin approximations a value very close to two is observed, while the $ECO$ determined for our Petrov-Galerkin approximations upgrades for $u_h$ and $\nabla \cdot {\bf p}_h$ and downgrades for ${\bf p}_h$. Since this significantly contrasts with the observations for Test-problems 4 and 5, we infer that, as a rule,  the performance of \eqref{Varimixh1} is quite acceptable and significantly superior to the one of its Galerkin counterpart in terms of accuracy.\\
Finally, this test-problem also shows that the polyhedral approximation of $\partial \Omega_{1/2}$ induces no significant erosion 
of both methods in terms of neither order nor accuracy. We recall that for their two-dimensional analogs studied in \cite{ESAIM-M2AN}, such an approximation of a boundary portion on which $u$ vanishes is supposed to introduce additional errors in the natural norms of the second order for $RT_k$, with $k \geq 1$, provided this function is sufficiently smooth. From the results in Table 8, one may presume that this  property remains true in the three-dimensional case.   
            
\begin{table*}[h!]
{\small 
\centering
\begin{tabular}{ccccccc} &\\ [-.3cm]  
$h$ & $\longrightarrow$ & $1/2$ & $1/4$ & $1/8$ & $1/16$ & $EOC$ 
\tabularnewline &\\ [-.3cm] \hline &\\ [-.3cm]
$\parallel {\bf p}_h - {\bf p} \parallel_h$     & $\longrightarrow$ & 0.85904684E-2 & 0.24696251E-2 & 0.67148753E-3 & 0.18855030E-3 & 1.8408
\tabularnewline &\\ [-.3cm] \hline &\\ [-.3cm] 
$\parallel \bar{\bf p}_h - {\bf p} \parallel_h$ & $\longrightarrow$ & 0.13098590E+0 & 0.32988206E-1 & 0.82619724E-2 & 0.20678983E-2 & 1.9953
\tabularnewline &\\ [-.3cm] \hline &\\ [-.3cm] 
\tabularnewline &\\ [-.3cm] \hline &\\ [-.3cm]
$\parallel \nabla \cdot({\bf p}_h - {\bf p}) \parallel_h$ & $\longrightarrow$ & 0.39055719E-1 & 0.10549432E-1 & 0.26887779E-2 & 0.67552705E-3 & 1.9532
\tabularnewline &\\ [-.3cm] \hline &\\ [-.3cm]  
$\parallel \nabla \cdot(\bar{\bf p}_h - {\bf p}) \parallel_h$ & $\longrightarrow$ & 0.72604599E-1 & 0.18275309E-1 & 0.45808071E-2 & 
0.11458562E-2 & 1.9953
\tabularnewline &\\ [-.3cm] \hline &\\ [-.3cm]
\tabularnewline &\\ [-.3cm] \hline &\\ [-.3cm] 
 $\parallel u_h - u \parallel_h$ & $\longrightarrow$ & 0.15615854E-1 & 0.43753666E-2 & 0.11239484E-2 & 0.28307537E-3 & 1.9318
\tabularnewline &\\ [-.3cm] \hline &\\ [-.3cm] 
$\parallel \bar{u}_h - u \parallel_h$  & $\longrightarrow$ & 0.63165871E-1 & 0.15551180E-1 & 0.38752427E-2 & 0.96787469E-3 & 2.0089 
\tabularnewline &\\ [-.3cm] \hline &\\ [-.3cm] 
\tabularnewline &\\ [-.3cm] \hline &\\ [-.3cm]
$N_{it}$ for $tol = 10^{-5}$ & $\longrightarrow$ & 10 & 7 & 6 & 25 & 
\tabularnewline &\\ [-.3cm] \hline &\\ [-.3 cm] 
\end{tabular}
%\vskip.2cm 
\caption{Errors for $RT_1$-solutions of Test-problem 6 \& number of iterations to attain the solution of \eqref{Varimixh1}} 
}
\label{table6}
\end{table*}

\end{enumerate}

\section{Stability considerations}

%This observation extends neither to problem \eqref{Varimixh} nor to problem \eqref{HRT0} 
We know that for any $\nu \geq 0$ the uniform stability is guaranteed for the Galerkin counterparts of both \eqref{Varimixh} (for any $k \geq 0$) and \eqref{HRT0}, according to well known results (cf. \cite{BrezziFortin} and \cite{JCAM}). However, as far as the Petrov-Galerkin formulation inherent to \eqref{Varimixh} and \eqref{HRT0} is concerned, it is not possible to be as assertive, since both problems are non symmetric, among other reasons. Nevertheless, on the basis of the analysis carried out in \cite{ESAIM-M2AN} for the two-dimensional case, it is reasonable to conjecture that for any $\nu \geq 0$ \eqref{Varimixh} is uniformly stable. Similarly this property should also hold for \eqref{HRT0}. As a matter of fact, the results displayed in Tables 1 through 5 corroborate such a conjecture, since the $EOC$s given therein reveal no significant erosion in the predictable orders of convergence, except perhaps for the order discrepancy of ca. $4$\%  in the case of the approximation of the divergence of ${\bf p}$, in Test-problem 3 only. However, we should remind in this connection that the degree of reliability of the $EOC$s is limited, owing to the relative coarseness of the meshes employed to compute them.\\
\indent On the other hand, provided $\nu >0$, the uniform stability of the numerical scheme associated with the Galerkin counterpart of \eqref{Varimixh1} is guaranteed, since in this case we are dealing with a coercive problem. Unfortunately this property in not so helpful to establish the uniform stability of problem \eqref{Varimixh1}, on the basis of the observation that it can be viewed as a small perturbation of its Galerkin counterpart. We recall that this fact was the main ingredient employed in previous author's work on Petrov-Galekin formulations to handle DOFs prescribed on curvilinear boundaries, such as \cite{ESAIM-M2AN}. The drawback here comes from the additional term 
$2 \nu^{-1}(\nabla \cdot {\bf r},\nabla \cdot \bar{\bf r})_h$ for a given ${\bf r} \in {\bf P}_h^1$ in the stability analysis, where $\bar{\bf r} \in {\bf Q}_h^1$ is the field associated with ${\bf r}$ in the way described in Subsection 4.2. Indeed, it seems very difficult, if not impossible, to bound $\|\nabla \cdot [{\bf r} - \bar{\bf r}] \|_h$ by the norm of ${\bf r}$ in ${\bf H}(div,\Omega_h)$ multiplied by an $O(h^{\alpha})$-term for a certain $\alpha >0$, as required (cf. \cite{ESAIM-M2AN}). Actually the downgrade of the $EOC$s shown in Tables 6 through 8, as compared to the expected order two for the $RT_1$-method, could be explained by a stability constant depending on $h$ like an $O(h^{\gamma})$-term for a small strictly positive $\gamma$. However we recall that \eqref{Varimixh1} is just a trick aimed at enabling less expensive linear system solving. Otherwise stated, the true problem to solve in practical applications would be of the same nature as \eqref{Varimixh}, eventually for $\nu =0$, rather than of the form \eqref{Varimixh1}. It turns out that the uniform stability of \eqref{Varimixh} for $k=1$ is expected to be a consequence of the uniform stability of its Galerkin counterpart, at least for sufficiently small $h$, as it is legitimate to conjecture from the arguments in \cite{ESAIM-M2AN}.\\
Anyhow, since it is not our purpose here to address theoretical issues, but rather show numerical evidence of our new approach's adequacy to improve accuracy, we endeavored to check its stability in computational terms, by solving \eqref{Varimixh1} for problems having the same solution $u$ with $f = \nu u - \Delta u$, letting the value of $\nu$ vary between one and $10^{-j}$ for quite a large positive integer $j$.\\
%in the natural norm of ${\bf H}(div\Omega_h) \times L^2(\Omega_h)$  

More precisely, the experiments carried out in this connection consist of solving a test-problem in the domain $\Omega_{hol}$ contained in   $\Omega_1$, which is either the unit ball $S$ or the ellipsoid $E$ with $a=0.6$, $b=0.8$ and $c=1$. We consider the exact solution $u=\eta^2  - \eta$ with $\Gamma_1=\partial \Omega_{1/2}$. In all cases we take $\nu$ successively equal to $1$, $10^{-2}$, $10^{-4}$, $10^{-6}$ and $10^{-8}$. The meshes for the sub-domain $\Omega^{+++}$ are constructed in the same way as for the previous test-problems, taking $L=2,4,8$ for the case of $S$ and $L=10$ for the case of $E$. In Tables 9, 10 and 11 we display respectively the errors of the approximations of ${\bf p}$, $\nabla \cdot {\bf p}$ and $u$ as solutions of \eqref{Varimix} in $\Omega_{hol}$ with $h=1/2$, $h=1/4$ and $h=1/8$ in the case of the unit ball $S$ and for $h=1/10$ in the case of the ellipsoid $E$.\\
The same iterative procedure described at the beginning of Subsection 4.2 was used to compute the approximate solutions. Except for the case of $\nu = 10^{-8}$ with $h > 1/2$, the tolerance of $10^{-5}$ for the increments of the normal flux DOFs of ${\bf p}_h$ was satisfied within less than $100$ iterations. \\
As one can infer from the three tables, the error behavior for a given pair mesh-domain is very stable as $\nu$ varies. Notice that this is true of all cases, including those where the iteration-stop criterion was not satisfied (except for a slightly larger discrepancy observed for $h=1/2$ and $\nu = 10^{-8}$; cf. Remark 3 below). This remarkable stability is worth being emphasized, since the importance in a theoretical stability analysis of the problematic term $2(\nabla \cdot {\bf p}_h,\nabla \cdot {\bf q})_h/\nu$ in formulation \eqref{Varimixh1} increases sharply as $\nu$ decreases. \\
Notice that, idealistically, a less unfair numerical stability check should be conducted in the more favorable case of \eqref{Varimixh} with values of $\nu$ all the way down to zero. However, should the assessment cost of such a formulation remain reasonable for $k=1$, the finest mesh in use would have to be even coarser than the one considered in Subsection 4.2. In this case the experiments might not be conclusive at all. 

\begin{table*}[h!]
{\small 
\centering
\begin{tabular}{ccccccccc} &\\ [-.3cm]  
$h \! \downarrow$ & \hspace{-2mm} $\Omega_1 \! \downarrow$ & \hspace{-2mm} $\nu$ \hspace{-5mm} & $\longrightarrow$ & $1$ & $10^{-2}$ & $10^{-4}$ & $10^{-6}$ & $10^{-8}$ 
\tabularnewline &\\ [-.3cm] \hline &\\ [-.3cm]
\tabularnewline &\\ [-.3cm] \hline &\\ [-.3cm]
$1/2$ & $S$ & & $\longrightarrow$  & 0.15565580E-1 & 0.15181653E-1 & 0.15185684E-1 & 0.15314842E-1 & 0.16593588E-1$^{+}$
\tabularnewline &\\ [-.3cm] \hline &\\ [-.3cm]
$1/4$ & $S$ & & $\longrightarrow$  & 0.46899481E-2 & 0.46080535E-2 & 0.46076749E-2 & 0.46076923E-2 & 0.46075754E-2$^{*}$
\tabularnewline &\\ [-.3cm] \hline &\\ [-.3cm]  
$1/8$ & $S$ & & $\longrightarrow$  & 0.14666198E-2 & 0.14522572E-2 & 0.14521662E-2 & 0.14521611E-2 & 0.14599640E-2$^{*}$
\tabularnewline &\\ [-.3cm] \hline &\\ [-.3cm]
\tabularnewline &\\ [-.3cm] \hline &\\ [-.3cm] 
$1/10$& $E$ & & $\longrightarrow$  & 0.12349943E-2 & 0.12321829E-2 & 0.12321913E-2 & 0.12321264E-2 & 0.12894992E-2$^{*}$
\tabularnewline &\\ [-.3cm] \hline &\\ [-.3cm] 
%\tabularnewline &\\ [-.3cm] \hline &\\ [-.3cm] 
\end{tabular}
%\vskip.2cm 
\caption{Errors $\|{\bf p}_h - {\bf p} \|_h$ for test-problems solved by \eqref{Varimixh1} with decreasing values of $\nu$} 
}
\label{table9}
\end{table*}
\begin{remark}
The errors $\|{\bf p}_h - {\bf p} \|_h$ in the cases where our iteration-stop criterion was not satisfied are marked with a $*$ on the right, while the error a little discrepant with respect to the one computed for other values of $\nu$ is marked with the sign $+$ on the right. \rule{2mm}{2mm}
\end{remark}

\begin{table*}[h!]
{\small 
\centering
\begin{tabular}{ccccccccc} &\\ [-.3cm]  
$h \! \downarrow$ & \hspace{-2mm} $\Omega_1 \! \downarrow$ & \hspace{-2mm} $\nu$ \hspace{-5mm} & $\longrightarrow$ & $1$ & $10^{-2}$ & $10^{-4}$ & $10^{-6}$ & $10^{-8}$ 
\tabularnewline &\\ [-.3cm] \hline &\\ [-.3cm]
\tabularnewline &\\ [-.3cm] \hline &\\ [-.3cm]
$1/2$ & $S$ & & $\longrightarrow$  & 0.45312145E-1 & 0.43476854E-1 & 0.43476666E-1 & 0.43477180E-1 & 0.43476663E-1
\tabularnewline &\\ [-.3cm] \hline &\\ [-.3cm]
$1/4$ & $S$ & & $\longrightarrow$  & 0.11388780E-1 & 0.10882607E-1 & 0.10882544E-1 & 0.10882544E-1 & 0.10882535E-1
\tabularnewline &\\ [-.3cm] \hline &\\ [-.3cm]  
$1/8$ & $S$ & & $\longrightarrow$  & 0.28575701E-2 & 0.27275941E-2 & 0.27275791E-2 & 0.27275806E-2 & 0.27275950E-2
\tabularnewline &\\ [-.3cm] \hline &\\ [-.3cm]
\tabularnewline &\\ [-.3cm] \hline &\\ [-.3cm] 
$1/10$& $E$ & & $\longrightarrow$  & 0.25416234E-2 & 0.25033932E-2 & 0.25033868E-2 & 0.25033869E-2 & 0.25034614E-2
\tabularnewline &\\ [-.3cm] \hline &\\ [-.3cm] 
%\tabularnewline &\\ [-.3cm] \hline &\\ [-.3cm] 
\end{tabular}
%\vskip.2cm 
\caption{Errors $\|\nabla \cdot ( {\bf p}_h - {\bf p}) \|_h$ for test-problems solved by \eqref{Varimixh1} with decreasing values of $\nu$} 
}
\label{table10}
\end{table*}

\begin{table*}[h!]
{\small 
\centering
\begin{tabular}{ccccccccc} &\\ [-.3cm]  
$h \! \downarrow$ & \hspace{-2mm} $\Omega_1 \! \downarrow$ & \hspace{-2mm} $\nu$ \hspace{-5mm} & $\longrightarrow$ & $1$ & $10^{-2}$ & $10^{-4}$ & $10^{-6}$ & $10^{-8}$ 
\tabularnewline &\\ [-.3cm] \hline &\\ [-.3cm]
\tabularnewline &\\ [-.3cm] \hline &\\ [-.3cm]
$1/2$ & $S$ & & $\longrightarrow$  & 0.15321424E-1 & 0.15745039E-1 & 0.15749417E-1 & 0.15761373E-1 & 0.15973668E-1
\tabularnewline &\\ [-.3cm] \hline &\\ [-.3cm]
$1/4$ & $S$ & & $\longrightarrow$  & 0.39944226E-2 & 0.41120834E-2 & 0.41133545E-2 & 0.41134668E-2 & 0.41131233E-2
\tabularnewline &\\ [-.3cm] \hline &\\ [-.3cm]
$1/8$ & $S$ & & $\longrightarrow$  & 0.10112449E-2 & 0.10417966E-2 & 0.10421020E-2 & 0.10421596E-2 & 0.10423046E-2
\tabularnewline &\\ [-.3cm] \hline &\\ [-.3cm]
\tabularnewline &\\ [-.3cm] \hline &\\ [-.3cm]   
$1/10$& $E$ & & $\longrightarrow$  & 0.50410167E-3 & 0.51347815E-3 & 0.51325951E-3 & 0.51402615E-3 & 0.51563509E-3
\tabularnewline &\\ [-.3cm] \hline &\\ [-.3cm] 
\end{tabular}
%\vskip.2cm 
\caption{Errors $\| u_h - u \|_h$ for test-problems solved by \eqref{Varimixh1} with decreasing values of $\nu$} 
}
\label{table11}
\end{table*}

\section{Concluding Remarks}
   
A rather large and varied amount of numerical experiments with the variants \eqref{Varimixh}, \eqref{Varimixh1} and \eqref{HRT0} of the $RT_k$ mixed method for $k=0$ and $k=1$, together with a Hermite analog of the $RT_0$-method, were reported in this article. The results were not fully conclusive as for the uniform stability of these methods. Additionally, superconvergence was observed in some cases. 
As far as we can see, any attempt to provide reasonable explanations for such unexpected or odd behaviors - even if favorable -, would be inhibited by basically two factors: the relative coarseness of the meshes employed in the experiments and the use of the enriched formulation \eqref{Varimixh1} in the case of the $RT_1$-method.  \\

Anyway, in spite of a few important open questions to be clarified in future work, we can assert that the main goal of the Raviart-Thomas method's variants introduced in this work was undoubtedly attained. Indeed, in most cases we have examined, the underlying modification designed to enforce prescribed normal fluxes on a curvilinear boundary portion, brought about a substantial enhancement of accuracy, as compared to the do-nothing approach. Furthermore, essentially in no case the latter showed up more accurate than the former. \\

More largely, on the basis of the results supplied in Section 4, as much as in previous author's work on the topic published in peer-review journals such as \cite{JAMP}, \cite{ZAMM}, \cite{IMAJNA}, \cite{JCAMbis} and \cite{ESAIM-M2AN}, besides other means, there is at least one conclusion to be taken for granted about the type of formulation advocated here: \\%  to handle DOFs prescribed on curvilinear boundaries: 
\indent It showed once again to be a simple and accurate alternative at a time, not only to shifting essential conditions of different types prescribed on the boundary of a curved problem-definition domain $\Omega$, to the boundary of approximating polytopes $\Omega_h$ \footnote{Generally speaking, an approximating polytope is the union of polygonal or polyhedral elements in meshes fitting $\Omega$.}, but also to widespread techniques such as parametric elements.\\

To finalize the author would like to underline the universality of the principle the technique adopted in this article is based upon, to handle essential conditions prescribed on curvilinear boundaries. In this sense, its application is not restricted at all to finite element methods. Actually, successful adaptions thereof to other numerical methods commonly used in Engineering and Science, such as the discontinuous Galerkin method, can already be found in the literature (see e.g. \cite{DG2D} and references therein).      

\section*{Acknowledgement} 
\vspace*{-3mm} 
The author is thankful to Fleurianne Bertrand for helpful discussions.

\end{document}